\documentclass[twoside,10pt,leqno]{amsart}
\usepackage{amsfonts}
\usepackage{amsmath}
\usepackage{amsrefs}
\usepackage{amsthm}
\everymath={\displaystyle}

\begin{document}
\newtheorem{theorem}{Theorem}
\newtheorem{lemma}{Lemma}
\newtheorem{corollary}{Corollary}
\newtheorem{conjecture}{Conjecture}
\numberwithin{equation}{section}
\renewcommand{\thefootnote}{\fnsymbol{footnote}}
\newcommand{\dif}{\mathrm{d}}
\newcommand{\intz}{\mathbb{Z}}
\newcommand{\ratq}{\mathbb{Q}}
\newcommand{\natn}{\mathbb{N}}
\newcommand{\comc}{\mathbb{C}}
\newcommand{\rear}{\mathbb{R}}
\newcommand{\prip}{\mathbb{P}}
\newcommand{\uph}{\mathbb{H}}

\title{\bf Oscillations of Hecke Eigenvalues at Shifted Primes}
\author{Liangyi Zhao}
\date{\today}
\maketitle

\begin{abstract}
In this paper, we are interested in exploring the cancellation of Hecke eigenvalues twisted with an exponential sums whose amplitude is $\sqrt{n}$ at prime arguments.
\end{abstract}

\section{History and Introduction}

In this paper, we are interested in estimating an exponential sum over primes with square root amplitude twisted with Hecke eigenvalues.  More precisely, we want to have an estimate for the following sum
\begin{equation} \label{sumprim2}
S(N) \stackrel{\mathrm{def}}{=} \sum_{n \leq N} \lambda (n) \Lambda (n) e ( \alpha \sqrt{n} ), \; \mbox{with} \; \alpha >0.
\end{equation}
$\lambda (n)$ henceforth denotes the normalized Fourier coefficients of a cusp form $f(z)$ of weight $k \geq 12$ for the full modular group, and $f(z)$ is an eigenform of all the Hecke operators.  {\it Id est},
\begin{equation*}
f(z) = (cz+d)^k f \left( \frac{az+b}{cz+d} \right) \; \mbox{for all} \; \left( \begin{array}{cc} a & b \\ c & d \end{array} \right) \in SL_2(\intz),
\end{equation*}
and $T_n f = \lambda(n) f$ for all $n \geq 1$, where $T_n$ is the $n$-th Hecke operator.  These Hecke eigenvalues, $\lambda(n)$'s, agree with the coefficients in the Fourier series expansion of $f(z)$
\begin{equation*}
f(z) = \sum_{n=1}^{\infty} \lambda(n) n^{\frac{k-1}{2}} e(nz).
\end{equation*}
$\Lambda (n)$, as usual, denotes the von Mangoldt function and $\alpha>0$.  References on the subject of Hecke eigenvalues are abundant.  See \cite{HI3}, \cite{JPS} and \cite{JS2}. \newline

Estimation of $S(N)$ is of interest from two points of view.  First the sum
\begin{equation*}
\sum_{n \leq N} \Lambda (n) e ( \alpha \sqrt{n} )
\end{equation*}
has been an object of interest ever since the method of I. M. Vinogradov was first developed.  It was Vinogradov himself \cite{IMV} who showed the afore-mentioned sum is $O( N^{\frac{7}{8}+ \epsilon} )$ with the implied constant depending on $\alpha$ and $\epsilon$. \newline

Second, the size and oscillations of the Hecke eigenvalues themselves are objects of great interest.  By Rankin-Selberg method, one achieves the asymptotics
\begin{equation*}
\sum_{n \leq N} |\lambda(n)|^2 \sim cN,
\end{equation*}
as $N$ tends to $\infty$ and $c$ here is a positive constant that depends on $f(z)$.  Of course, we also have the following.
\begin{theorem} [Ramanujan Conjecture] \label{ramconj}
With $\lambda(n)$ denoting the $n$-th Hecke eigenvalue of a cusp form, $f(z)$, for the full modular group, we have
\begin{equation*} 
|\lambda(n)| \leq \tau(n) \ll n^{\epsilon},
\end{equation*}
where $\tau(n)$ is the divisor function.
\end{theorem}

\begin{proof} This famous result was of course, proved by P. Deligne \cite{PD1} in 1974, and we shall appeal to this theorem in later sections. \end{proof}

Regarding the sign changes of the Hecke eigenvalues, it was due to Hardy and Ramanujan, and A. Good \cite{AG}, respectively that
\begin{equation*}
\sum_{n \leq N} \lambda (n) e ( \alpha n ) \ll N^{\frac{1}{2}} \log (2N), \; \mbox{and} \sum_{n \leq N} \lambda (n) \ll N^{\frac{1}{3}+\epsilon}.
\end{equation*}
Moreover, M. R. Murty \cite{MRM} conjectured $\Omega$ result that
\begin{equation*}
\sum_{p \leq N} \lambda (p) = \Omega_{\pm} \left( \frac{\sqrt{N} \log
\log \log N}{\log N} \right).
\end{equation*}
and succeeded in proving it provided some $L$-function has no real zero between 1/2 and 1.  Also, S. D. Adhikari \cite{SDA} proved essentially the same result for cusp forms for the group $\Gamma_0(N)$. \newline

The method that we employ in estimating $S(N)$ is that developed by Vinogradov \cite{IMV}.  As one familiar with the
method knows, the best possible results that technology can prove is $S(N) =O\left( N^{\Theta + \epsilon} \right)$, with $\Theta = \frac{3}{4}$.  However, by the so-called ``principle of square-rooting,'' then one may be led to believe that $S(N) =O\left( N^{\frac{1}{2} + \epsilon} \right)$.  Surprisingly, Iwaniec, Luo and Sarnak \cite{ILS} gave conditional heuristics that $\Theta = \frac{3}{4}$ is where truth actually lies.  They have the asymptotic formula, under the assumption of some extremely strong hypotheses,
\begin{equation} \label{ilsheur}
S(N) = Z N^{\frac{3}{4}} + O(N^{\frac{5}{8}+\epsilon}),
\end{equation}
where $Z$ is a non-zero constant that depends on the cusp form.  In this paper, we prove the following.

\begin{theorem} \label{sumprimresult}
With $S(N)$ defined as in \eqref{sumprim2}, we have
\begin{equation*}
S(N) \ll N^{\frac{5}{6}} [\log (3N)]^{21},
\end{equation*}
where the implied constant depends effectively on $\alpha$ in \eqref{sumprim2} and the cusp form $f(z)$.
\end{theorem}

The author wishes to thank his thesis advisor, Professor Henryk Iwaniec, for suggesting this problem and who, with his advise and support, has been most generous.  Moreover, the author also thanks the referee for his many helpful comments on the original manuscript. \newline

The following notations and conventions are used throughout paper. \newline

\noindent $e(z) = \exp (2 \pi i z) = e^{2 \pi i z}$. \newline
$f = O(g)$ means $|f| \leq cg$ for some unspecified postive constant $c$. \newline
$f \ll g$ means $f=O(g)$.  Unless otherwise stated, all implied constants in $\ll$ and $O$ are absolute. \newline
$L(s,f) = \sum_{n=1}^{\infty} \lambda(n) n^{-s}$ is the automorphic $L$ function for a cusp form $f$, with $\lambda(n)$ being the same as those in \eqref{sumprim2}. \newline
$\qed$ denotes the end of a proof or the proof is easy and standard.

\section{Preliminary Lemmas}

In this section, we quote the results needed later.  First, we have the multiplicative property of Hecke eigenvalues.

\begin{lemma} \label{multhecke}
Hecke eigenvalues are multiplicative and they satisfy the
following relation.
\begin{equation*}
\lambda (mn) = \sum_{d| \gcd(m,n)} \mu (d) \lambda \left( \frac{m}{d} \right) \lambda \left( \frac{n}{d} \right).
\end{equation*}
\end{lemma}

\begin{proof} This Lemma follows by applying the M\"{o}bius inversion formula to the product formula for the Hecke eigenvalues.  See, for example, Proposition 14.9 of \cite{HIEK}. \end{proof}

We shall also need to estimate the average of the divisor functions.

\begin{lemma} \label{divisorave2}
With $k, l \in \natn$ and with $\tau_k(n)$ denoting the number of ways $n$ can be written as product of $k$ integers, then
\begin{equation*}
\sum_{n \leq x } [\tau_k(n) ]^l \ll x [\log (2x) ]^{k^l-1},
\end{equation*}
where the implied constant depends on $k$.
\end{lemma}
\begin{proof} The proof is easy and standard.  See (1.80) of \cite{HIEK}.  \end{proof}

We shall need the mean value theorem of the square of automorphic $L$-functions on the critical line.

\begin{lemma}
We have
\begin{equation*}
\int_0^T \left| L \left( \frac{1}{2}+ it, f \right) \right|^2 \dif t \ll T^{1+\epsilon},
\end{equation*}
where the implied constant depends on $\epsilon$ alone.
\end{lemma}
\begin{proof} The result arrives via similar means as the analogous result for the Riemann Zeta-function.  One can, of course, prove stronger results, as in \cite{AG2}, but the above suffices.  Similar result also holds for $L' \left( \frac{1}{2} + it, f \right)$. \end{proof}

In our proof, we shall need to estimate certain exponential sums.  The following lemmas suffice for our enterprise. \newline

\begin{lemma} \label{specpoiss}
Let $f(x)$ be a real-valued function with $|f'(x)| \leq 1- \theta$ and $f''(x) \neq 0$ on $[a,b]$.  We have
\begin{equation*}
\sum_{a <n<b} e [f(n)] = \int_a^b e[f(x)] \dif x +O( \theta^{-1}),
\end{equation*}
where the implied constant is absolute.
\end{lemma}

\begin{proof} This is a special case of the Truncated Poisson Summation Formula in \cite{ET} and \cite{HM3}. \end{proof}

Next we have these estimates for exponential integrals.
\begin{lemma} \label{Bproc1}
Let $f(x)$ be a real-valued function with two continuous derivatives on $[a,b]$ and that $f'(x)f''(x) \neq 0$ on $[a,b]$, then we have
\begin{equation*}
\left| \int_a^b e[f(x)] \dif x \right| \leq \pi^{-1} \left( |f'(a)|+|f'(b)| \right).
\end{equation*}
\end{lemma}
\begin{proof} This result arrives easily from partial integration.  See Lemma 8.9 of \cite{HIEK}. \end{proof}

We need a result known as stationary phase which we get to through
\begin{lemma} \label{statphase1}
Let $h(x)$ be a real function with two continuous derivatives on $[0,X]$ such that
\begin{equation} \label{statphase2}
h(0)=1, \; h(x) \gg 1, \; \left[ x h(x) \right]' \gg 1, \; h'(x) \ll \frac{1}{X}, \; \mbox{and} \; h''(x) \ll \frac{1}{X^2}.
\end{equation}
Then for $\alpha > 0$, we have
\begin{equation} \label{statphase3}
\int_0^X e[\alpha x^2 h(x) ] \dif x = \frac{8}{\sqrt{\alpha}} e \left( \frac{1}{8} \right) + O \left( \frac{1}{\alpha X} \right),
\end{equation}
where the implied constant in \eqref{statphase3} depends on those in \eqref{statphase2}.
\end{lemma}
\begin{proof} This lemma is proved by standard means.  See Lemma 8.14 of \cite{HIEK}. \end{proof}

With the previous lemma at our disposal, we have the following.

\begin{lemma} [Stationary Phase] \label{vandercor2}
Let $f(x)$ be a real valued function with four continuous derivatives on $[a,b]$ with
\begin{equation*}
f''(x) \geq \Lambda, \; |f'''(x)| \leq \Lambda X^{-1}, \; \mbox{and} \; |f^{(4)}(x)| \leq \Lambda X^{-2},
\end{equation*}
for some $\Lambda >0$ and $X>0$.  Also suppose $f'(c)=0$ with $c \in (a,b)$.  Then
\begin{equation*}
\int_a^b e[f(x)] \dif x = e \left[ f(c)+ \frac{1}{8} \right] f''(c)^{-\frac{1}{2}} + O \left[ \frac{1}{\Lambda} \left( \frac{1}{b-c}+\frac{1}{c-a}+\frac{1}{X} \right) \right], 
\end{equation*}
where the implied constant is absolute.
\end{lemma}
\begin{proof} This lemma follows by Lemma~\ref{statphase1} and the second degree Taylor approximation of $f(x)$.  This is Corollary 8.15 in \cite{HIEK}. \end{proof}

We shall also need the following Perron-type formula which approximates Dirichlet polynomials.

\begin{lemma} [Perron's Formula] \label{perronform}
Let $f(s)$ be defined by the Dirichlet series
\[ f(s) = \sum_{n=1}^{\infty} \frac{a_n}{n^s}, \; \mbox{for} \; \Re s >1. \]
where $a_n \ll \psi(n)$ for some non-decreasing function $\psi(n)$, and
\[ \sum_{n=1}^{\infty} \frac{a_n}{n^s} = O \left[ \frac{1}{(\sigma-1)^{\alpha}} \right], \]
as $\sigma$ tends to $1$.  Then if $c >0$, $\sigma +c>1$, and $x \notin \natn$, we have
\begin{equation*}
\sum_{n <x} \frac{a_n}{n^s} = \frac{1}{2 \pi i} \int_{c-iT}^{c+iT} f(s+w) \frac{x^w}{w} \dif w + O \left[ \frac{\psi(x) x^{1-\sigma}}{T \| x \|} \frac{x^c}{T(\sigma+c-1)^{\alpha}} + \frac{\psi(2x) x^{1-\sigma} \log x }{T} \right];
\end{equation*}
where $\| x \| = \min_{k \in \intz} |x-k|$, and if $x \in \natn$, then
\begin{equation*}
\sum_{n=1}^{x-1} \frac{a_n}{n^s} + \frac{a_x}{2x^s} = \frac{1}{2 \pi i} \int_{c-iT}^{c+iT} f(s+w) \frac{x^w}{w} \dif w + O \left[ \frac{x^c}{T(\sigma+c-1)^{\alpha}} + \frac{\psi(2x) x^{1-\sigma} \log x}{T} + \frac{\psi(x) x^{-\sigma}}{T} \right].
\end{equation*}
\end{lemma}
\begin{proof} This is quoted from \cite{ET}. \end{proof}

Finally, we need an integral form of the large sieve inequality in the estimation of the mean values of Dirichlet polynomials.

\begin{theorem} [Large Sieve] \label{intlarsiev}
Suppose that $\lambda_1, \cdots \lambda_N$ are distinct real numbers, and supposed that $\delta >0$ is chosen so that $|\lambda_m - \lambda_n| \geq \delta, \; \mbox{for} \; m \neq n$.  Then for any complex coefficients $a_n$, and any $T>0$, we have
\begin{equation*}
\int_0^T \left| \sum_{n=1}^N a_n e (\lambda_n t ) \right|^2 \dif t = \left( T + \frac{\theta}{\delta} \right) \sum_{n=1}^N |a_n|^2
\end{equation*}
for some $\theta$ with $ -1 \leq \theta \leq 1$.
\end{theorem}

\begin{proof} This is quoted from \cite{HM3} and derived using Selberg's majorant and minorant. \end{proof}

Theorem~\ref{intlarsiev} is applicable to Dirichlet polynomials by the mean value theorem of differential calculus.  We have
\begin{equation*}
\int_0^T \left| \sum_{n=1}^N a_n n^{-it} \right|^2 \dif t = \left[ T + O(N) \right] \sum_{n=1}^N |a_n|^2,
\end{equation*}
where the implied constant is absolute.

\section{Partition of the von Mangoldt Function}

We begin with the following identity.

\begin{lemma} [Vaughan]
Suppose $y \geq 2$ is a real positive number, then if $n>y$, we have
\begin{equation} \label{vaughanid}
\Lambda(n) = \sum_{\substack{ab=n \\ b \leq y}} \mu (b) \log a - \sum_{\substack{abc=n \\ b,c \leq y}} \mu (b) \Lambda(c) + \sum_{\substack{abc=n \\ b,c >y}} \mu (b) \Lambda(c).
\end{equation}
\end{lemma}
\begin{proof} This is the identity in \cite{RCV}. \end{proof}

We shall use a variant of the above idenity.  More precisely, we have
\begin{lemma} \label{vaughanid2}
For $N < n \leq 2N$, $n \geq y = N^{\frac{1}{3}}$ and $z = \sqrt{2N}$, we have
\begin{equation*}
\Lambda(n) = \Lambda_1(n) +\Lambda_2(n)+\Lambda_3(n) + \Lambda_4(n)+\Lambda_5(n),
\end{equation*}
where
\begin{equation*}
\Lambda_1(n) = \sum_{\substack{ab=n \\ b \leq y}} \mu (b) \log a, \; \Lambda_2(n) = -\sum_{\substack{abc=n \\ b,c \leq y \\ a \geq z}} \mu (b) \Lambda(c), \; \Lambda_3(n) = - \sum_{\substack{abc=n \\ b,c \leq y \\ y < a < z}} \mu (b) \Lambda(c),
\end{equation*}
\begin{equation*}
\Lambda_4(n) = \sum_{\substack{abc=n \\ c >y \\ y < b < z}} \mu (b) \Lambda(c), \; \mbox{and} \; \Lambda_5(n) = \sum_{\substack{abc=n \\ b \geq z \\ y < c \leq z}} \mu (b) \Lambda(c).
\end{equation*}
\end{lemma}
\begin{proof} We decompose the right-hand side of \eqref{vaughanid} further.  Take $N < n \leq 2N$ in dyadic intervals, $y = N^{\frac{1}{3}}$ and $z = \sqrt{2N}$.  We have
\begin{equation} \label{part1}
\sum_{\substack{abc=n \\ b,c \leq y}} \mu (b) \Lambda(c) = \sum_{\substack{abc=n \\ b,c \leq y \\ a \geq z}} \mu (b) \Lambda(c) + \sum_{\substack{abc=n \\ b,c \leq y \\ y < a < z}} \mu (b) \Lambda(c).
\end{equation}

The partition according to the dichotomy of either $a \geq z$ or $a < z$ is obvious.  The extra condition in the second term of the right-hand side of \eqref{part1} that $y<a$ is due to $b,c \leq y \Longrightarrow bc \leq y^2 = N^{\frac{2}{3}}$, together with $abc = n \leq 2N$, we have $ a > N^{\frac{1}{3}} =y$. \newline

The third sum in \eqref{vaughanid} is similarly decomposed.  We have
\begin{equation} \label{part2}
\sum_{\substack{abc=n \\ b,c >y}} \mu (b) \Lambda(c) = \sum_{\substack{abc=n \\ c >y \\ y < b < z}} \mu (b) \Lambda(c) + \sum_{\substack{abc=n \\ b \geq z \\ y < c \leq z}} \mu (b) \Lambda(c).
\end{equation}

Again, the dichotomy of $b <z$ or $b \geq z$ in the right-hand side of \eqref{part2} is obvious.  Furthermore, the extra condition in the second sum of the right-hand side of \eqref{part2} is apparent as
\[ abc =n < N \; \mbox{and} \; b \geq z = \sqrt{2N} \Longrightarrow ac \leq \sqrt{2N} = z \Longrightarrow c \leq z. \]

Therefore, combining Lemma~\ref{vaughanid} and \eqref{part1} and \eqref{part2}, we have the desired result. \end{proof}

Thus the sum of our interest in \eqref{sumprim2} is decomposed and it suffices to estimate each individual component.  We have
\begin{equation} \label{part4}
\left| \sum_{N < n \leq 2N} \Lambda(n) \lambda(n) e(\alpha \sqrt{n}) \right| \leq \sum_{i=1}^5 \left| S_i(N) \right|,
\end{equation}
where, $S_i(N) = \sum_{N < n \leq 2N} \Lambda_i(n) \lambda(n) e (\alpha \sqrt{n})$, and $\Lambda_i(n)$'s are defined in Lemma~\ref{vaughanid2}.

\section{Bilinear Forms Treatment}

The last three sums of \eqref{part4} are similar and can be disposed using similar means.  Toward that end, we have
\begin{lemma} \label{biline}
With $S_i(N)$ defined as before, we have
\begin{equation*}
|S_3(N)| + |S_4(N)| + |S_5(N)| \ll N^{\frac{5}{6}} (\log 3N)^{20},
\end{equation*}
where the implied constant depends on $\alpha$ and the cusp form $f(z)$.
\end{lemma}
\begin{proof} Breaking the summations into dyadic intervals, it suffices to estimate, for arithmetic functions $\beta(m)$ and $\gamma(l)$, sums of the following shape,
\begin{equation} \label{biline1}
\sum_{M < m \leq 2M} \beta(m) \left| \sum_{L < l \leq 2L} \gamma(l) \lambda (lm) e(\alpha \sqrt{lm}) \right|,
\end{equation}
with $N^{\frac{1}{2}} \leq M \leq N^{\frac{2}{3}}$, $N^{\frac{1}{3}} \leq L \leq N^{\frac{1}{2}}$, $ML =N$. \newline

Applying the multiplicative properties of Hecke eigenvalues, Lemma~\ref{multhecke}, \eqref{biline1} becomes
\begin{eqnarray*}
 & \sum_{M < m \leq 2M} \beta(m) \left| \sum_{s |m} \mu(s) \lambda \left( \frac{m}{s} \right) \sum_{\frac{L}{s} < t \leq \frac{2L}{s}} \gamma(st) \lambda (t) e(\alpha \sqrt{stm}) \right| \\
\leq & \sum_{M < m \leq 2M} \beta(m) \sum_{s |m} \left| \lambda \left( \frac{m}{s} \right) \right| \left| \sum_{\frac{L}{s} < t \leq \frac{2L}{s}} \gamma(st) \lambda (t) e(\alpha \sqrt{stm}) \right|.
\end{eqnarray*}

We divide the range of summation of the inner-most sum of the above further and estimate sums of the shape
\begin{equation} \label{beline11}
\sum_{M < m \leq 2M} \beta(m) \sum_{s |m} \left| \lambda \left( \frac{m}{s} \right) \right| \left| \sum_{\frac{L}{s} < t \leq \frac{L+L_0}{s}} \gamma(st) \lambda (t) e(\alpha \sqrt{stm}) \right|,
\end{equation}
where $L_0 \leq L$ will be chosen later.  We note that the number of sums like the above is $O \left[ \frac{L}{L_0} \left( \log N \right)^2 \right]$.  We apply the Cauchy-Schwarz inequality and majorize $\sum_{s|m} \left| \lambda \left( \frac{m}{s} \right) \right|^2$ by $\tau(m)^3$, we see that the square of the sum in \eqref{beline11} is majorized by
\begin{equation*}
\left( \sum_{M < m < 2M} \tau(m)^3 \left| \beta(m) \right|^2 \right) \left( \sum_{s \leq 2M} \sum_{\frac{M}{s} < m \leq \frac{2M}{s}} \left| \sum_{\frac{L}{s}< t \leq \frac{L+L_0}{s}} \gamma(st) \lambda(t) e(\alpha s \sqrt{tm}) \right|^2 \right).
\end{equation*}
Opening up the complex modulus square in the second factor and swapping the order of summations, it becomes
\begin{equation} \label{biline3}
\sum_{s \leq 2M} \mathop{\sum \sum}_{\frac{L}{s} < t, t' \leq \frac{L+L_0}{s}} \gamma(st) \overline{\gamma}(st') \lambda(st) \overline{\lambda}(st') \sum_{\frac{M}{s} < m \leq \frac{2M}{s}} e [ \alpha s (\sqrt{t}-\sqrt{t'}) \sqrt{m} ].
\end{equation}
The contribution of the diagonal terms in \eqref{biline3} is
\[ \ll M \sum_{s \leq 2M} s^{-1} \sum_{\frac{L}{s} < t \leq \frac{L+L_0}{s}} \left| \gamma(st) \lambda(st) \right|^2. \]
For the terms with $t \neq t'$, we note that if $f(m) = \alpha s (\sqrt{t}-\sqrt{t'})\sqrt{m}$, then 
\[ f'(m) = \frac{\alpha s (\sqrt{t}-\sqrt{t'})}{2\sqrt{m}} = \frac{\alpha s (t-t')}{2\sqrt{m}(\sqrt{t}+\sqrt{t'})} \leq \frac{\alpha s L_0}{2\sqrt{N}}. \]
Therefore, by choosing $L_0 = \frac{7L}{8 \alpha s}$ if $\alpha \geq 1$ or $L_0 = \frac{7L}{8s}$ if $0 < \alpha <1$ and recalling that $L_0 \leq L \leq \sqrt{N}$, we can ensure that $f'(m) \leq \frac{7}{16}$ for all $m$'s of interest.  The inner-most sum of \eqref{biline3} is well-approximated by that of its corresponding integral by Lemma~\ref{specpoiss}.  The modulus in question is
\begin{equation*}
= \int_{\frac{M}{s}}^{\frac{2M}{s}} e \left[ \alpha s (\sqrt{t}-\sqrt{t'})\sqrt{x} \right] \dif x + O(1) \ll \frac{\sqrt{M}}{s^{\frac{3}{2}} (\sqrt{t}-\sqrt{t'})} +1,
\end{equation*}
where the last implied constant depends on $\alpha$ and the last inequality arrives by the virtue of Lemma~\ref{Bproc1}.  Observe that $\beta(m)$ is
\begin{equation} \label{coef1}
 \sum_{\substack{bc=m \\ b,c\leq y}} \mu(b) \Lambda(c), \; \sum_{\substack{c|m \\ c >y}} \Lambda(c) \; \mbox{and} \; \sum_{\substack{b|m \\ b \geq z}} \mu(b)
\end{equation}
for $S_3(N)$, $S_4(N)$ and $S_5(N)$ respectively.  The moduli of three sums in \eqref{coef1} are majorized by $\log (2m)$, $\log (2m)$ and the divisor function $\tau (m)$, respectively.  In all three cases, we have,
\begin{equation*}
\sum_{M < m \leq 2M} \tau^3 (m) \left| \beta(m) \right|^2 \ll M [ \log (3M)]^{31},
\end{equation*}
with implied constant absolute, by the virtue of Lemma~\ref{divisorave2}. \newline

Similarly, $\gamma(l)$ is $1$, $\mu(l)$ and $\Lambda(l)$ for $S_3(N)$, $S_4(N)$ and $S_5(N)$ respectively.  In all three cases, $
\left| \gamma(l) \right| \leq \log (3l)$.  Consequently, \eqref{biline3} is majorized by
\begin{equation} \label{biline4}
\ll M^2 [\log (3N)]^{31} \sum_{s \leq 2M} s^{-1} \sum_{\frac{L}{s} < t \leq \frac{L+L_0}{s}} \left| \lambda(st) \right|^2 +M[\log (3N)]^{31} \sum_{s \leq 2M} \mathop{\sum \sum}_{\substack{\frac{L}{s} < t, t' \leq \frac{L+L_0}{s} \\ t \neq t'}} \left[ \frac{\sqrt{M}}{s^{\frac{3}{2}}| \sqrt{t}-\sqrt{t'} |} + 1 \right].
\end{equation}

Using the Ramanujan conjecture, Theorem~\ref{ramconj}, \eqref{biline4} is majorized by
\[ \ll \left[ M^2 L + M^{\frac{3}{2}} L^{\frac{3}{2}} + ML^2 \right] (\log 3N)^{35} \ll N^2 \left[ \frac{1}{L} + \frac{1}{N^{\frac{1}{2}}} + \frac{1}{M} \right] (\log 3N)^{35}. \]
Recall that $N^{\frac{1}{2}} \leq M \leq N^{\frac{2}{3}}$, $N^{\frac{1}{3}} \leq L \leq N^{\frac{1}{2}}$ and $ML=N$.  After taking the square root and then add up the sums over all the dyadic intervals, we have the desired result. \end{proof}

\section{Type I Sums}

It still remains to estimate the other terms in \eqref{part4} which are similar.  We take $f(t) = \frac{t}{2\pi} \log \left( \frac{t}{ex} \right)$.  $f(t)$ satisfies the conditions of Lemma~\ref{vandercor2} with $\Lambda= (2 \pi T)^{-1}$ and $c=x$.  Therefore, by the virtue of the said lemma,
\begin{equation*}
\int_T^{4T} e[f(t)] \dif t = e \left( \frac{-x}{2\pi} + \frac{1}{8} \right) \sqrt{2\pi x} + O\left[ T \left( \frac{1}{4T-x} + \frac{1}{x-T} + \frac{1}{T} \right) \right],
\end{equation*}
where $T=\pi \alpha \sqrt{N}$ and $x = 2\pi \alpha \sqrt{n}$, with $N < n \leq 2N$.  Note $2T \leq x \leq 3T$.  The above expression yields
\begin{equation} \label{trans3}
e( -\alpha \sqrt{n} ) = e \left( - \frac{1}{8} \right) (2\pi \sqrt{\alpha} )^{-1} n^{-\frac{1}{4}} \int_T^{4T} e[f(t)] \dif t + O \left( N^{-\frac{1}{4}} \right).
\end{equation}

We dispose of $S_1(N)$ with the following Lemma.

\begin{lemma} \label{S1est}
With $S_1(N)$ defined as before and for any $\epsilon >0$ given, we have
\begin{equation*}
S_1(N) \ll N^{\frac{3}{4}+\epsilon},
\end{equation*}
where the implied constant depends on $\alpha$, the cusp form $f(z)$ and $\epsilon$.
\end{lemma}

\begin{proof} By the virtue of \eqref{trans3}, it suffices to estimate
\begin{equation} \label{trans11}
e \left( - \frac{1}{8} \right) (2\pi \sqrt{\alpha} )^{-1} \int_T^{4T} \sum_{N <n \leq 2N} \Lambda_1(n) \lambda(n) n^{-\frac{1}{4}} \left( \frac{t}{2 \pi e \alpha \sqrt{n}} \right)^{it} \dif t + O \left[ N^{\frac{3}{4}} (\log 2N)^4 \right].
\end{equation}
After applying partial summation to integrand of \eqref{trans11}, it suffices to estimate, for $N < M \leq 2N$,
\begin{equation} \label{trans12}
N^{-\frac{1}{4}} \int_T^{4T} \left| \sum_{N < n \leq M} \Lambda_1(n) \lambda(n) n^{-\frac{it}{2}} \right| \dif t = 2 N^{-\frac{1}{4}} \int_{\frac{T}{2}}^{2T} \left| \sum_{N < n \leq M} \sum_{\substack{ab=n \\ b \leq y}} \mu(b) \log a \lambda(ab) (ab)^{-it} \right| \dif t.
\end{equation}

Applying Lemma~\ref{multhecke} and swapping the order of summations of the integrand of \eqref{trans12}, it becomes
\begin{equation} \label{trans13}
 \left| \sum_{d \leq y} \frac{\mu(d)}{d^{2it}} \sum_{l \leq \frac{y}{d}} \frac{\lambda(l) \mu(ld)}{l^{it}} \sum_{\frac{N}{ld^2} < h \leq \frac{M}{ld^2}} \frac{\lambda(h) (\log h + \log d)}{h^{it}} \right|.
\end{equation}
Therefore, the Dirichlet series we must consider is that of
\[ \left[ \log d L(s,f)-L'(s,f) \right] \sum_{l \leq \frac{y}{d}} \frac{\mu(ld) \lambda(l)}{l^s}, \]
where $L(s,f)$ is the $L$-function for the cusp form $f(z)$.  We first consider the sum with $L'(s,f)$ and the sum with $L(s,f)$ is treated similarly and yields the same majorant.  By the virtue of Perron's formula, Lemma~\ref{perronform}, we have
\begin{equation*}
\sum_{l \leq \frac{y}{d}} \frac{\mu(ld) \lambda(l)}{l^{it}} \sum_{\frac{N}{ld^2} < h \leq \frac{M}{ld^2}} \frac{\lambda(h) \log h}{h^{it}} = \frac{-1}{2\pi i} \int_{\sigma-iT}^{\sigma+iT} L'(w+it,f) \left[ \sum_{l \leq \frac{y}{d}} \frac{\mu(ld) \lambda(l)}{l^{w+it}} \right] \frac{M^w-N^w}{d^{2w}w} \dif w + O \left( \frac{N^{1+\epsilon}}{T} \right),
\end{equation*}
where, in the notations of Lemma~\ref{perronform}, we take $\sigma = 1 + \frac{1}{\log N}$, $s=it$ and $\psi (x) \ll x^{\epsilon}$ can be chosen.  Insert the above into \eqref{trans13}, the part corresponding to $L'(s,f)$ is 
\begin{equation*}
\frac{1}{2 \pi i} \sum_{d \leq y} \frac{\mu(d)}{d^{2it}} \int_{1 + \frac{1}{\log N}-iT}^{1 + \frac{1}{\log N}+iT} L'(w+it,f) \left( \sum_{l \leq \frac{y}{d}} \frac{\mu(ld) \lambda(l)}{l^{w+it}} \right) \frac{M^w-N^w}{d^{2w}w} \dif w + O \left( \frac{N^{1+\epsilon}}{T} \right).
\end{equation*}
Consider the contour given by the rectangle whose vertices are $\frac{1}{2} \pm iT$ and $1+\frac{1}{\log N} \pm iT$.  The above expression is well-approximated by
\begin{equation} \label{trans16}
\frac{1}{2 \pi i} \sum_{d \leq y} \frac{\mu(d)}{d^{2it}} \int_{\frac{1}{2}- iT}^{\frac{1}{2}+iT} L'(w+it,f) \left( \sum_{l \leq \frac{y}{d}} \frac{\mu(ld) \lambda(l)}{l^{w+it}} \right) \frac{M^w-N^w}{d^{2w}w} \dif w + O \left( \sqrt{\frac{N}{T}} y^{\frac{1}{2}+\epsilon} \right),
\end{equation}
where the second term in \eqref{trans16} is to majorize the contribution of the contour integral on the horizontal line segments of the rectangle and comes from convexity bound of $L'(s,f)$ on the critical line, $L' \left( \frac{1}{2}+it , f \right) \ll t^{\frac{1}{2}+\epsilon}$ and trivial bounds over the other factors.  Subconvexity bounds are known for these $L$-functions, see \cites{AG2, AG3}, but the trivial bound suffices for our purpose. \newline

Inserting everything into \eqref{trans12}, applying H\"{o}lder's inequality, $S_1(N)$ is majorized
\begin{equation*}
\ll N^{\frac{1}{4}} \sum_{d \leq y} \frac{1}{d} \int_{-iT}^{iT} \left[ \int_{\frac{T}{2}}^{2T} \left| L'(\frac{1}{2}+it,f) \right|^2 \dif t \right]^{\frac{1}{2}} \left[ \int_{\frac{T}{2}}^{2T} \left| \sum_{l \leq \frac{y}{d}} \frac{\mu(ld) \lambda(l)}{l^{\frac{1}{2}+it}} \right|^2 \dif t \right]^{\frac{1}{2}} \frac{\dif w}{1+2|w|} + N^{\frac{3}{4}+\epsilon}.
\end{equation*}
We apply Theorem~\ref{intlarsiev} to the second inner integral and it is $O\left[\left( T + \frac{y}{d} \right) \left( \frac{y}{d} \right)^{\epsilon} \right]$.  We utilize the mean value theorem for automorphic $L$-functions for the first inner integral and it is $O(T^{1+\epsilon})$.  Recall that $T = \pi \alpha \sqrt{N}$ and $y=N^{\frac{1}{3}}$.  The contribution in \eqref{trans13} of $L(s,f)$ is treated precisely the same way as the treatment of $L'(s,f)$.  Inserting those bounds, Lemma~\ref{S1est} follows. \end{proof}

To estimate the size of $S_2(N)$, we have
\begin{lemma} \label{S2est}
With $S_2(N)$ defined as before and for $\epsilon >0$ given, we have
\begin{equation*}
S_2(N) \ll N^{\frac{3}{4}+\epsilon},
\end{equation*}
where the implied constant depends on $\alpha$, the cusp form $f(z)$ and $\epsilon$.
\end{lemma}

\begin{proof} Consider a type of truncated M\"{o}bius convolution as follows.
\begin{equation*}
\nu (d) = \sum_{\substack{bc=d \\ b,c \leq y}} \mu(b) \Lambda(c).
\end{equation*}
It is elementary that $|\nu(d)| \leq \log 2d$.  We decompose $\Lambda_2(n)$ as follows.
\begin{equation*}
\Lambda_2(n) = - \sum_{\substack{ad=n \\ z > d}} \nu(d) + \sum_{\substack{ad=n \\ a,d < z}} \nu(d).
\end{equation*}

Therefore, we have also correspondingly decomposed $S_2(N)$ further and we have that $|S_2(N)|$ is bounded above by
\begin{equation*}
\left| \sum_{N < n \leq 2N} \sum_{\substack{ad=n \\ z > d}} \nu(d) \lambda(ad) e( \alpha \sqrt{ad} ) \right| + \left| \sum_{N < n \leq 2N} \sum_{\substack{ad=n \\ a,d < z}} \nu(d) \lambda(ad) e( \alpha \sqrt{ad} )  \right|.
\end{equation*}

For convenience, let the first sum in right-hand side of the above be $S_{21}(N)$ and the second $S_{22}(N)$.  Inserting \eqref{trans3} into the definition of $S_{21}(N)$ and applying partial summation to the integrand, we arrive at the following
\begin{equation} \label{trans211}
N^{-\frac{1}{4}} \int_{\frac{T}{2}}^{2T} \left| \sum_{N < n \leq M} \sum_{\substack{ad=n \\ z > d}} \nu(d) \lambda(ad) (ad)^{-it} \right| \dif t + N^{\frac{3}{4}+\epsilon}, \; \mbox{with} \; N < M \leq 2N.
\end{equation}
We now apply the multiplicative property of Hecke eigenvalues, Lemma~\ref{multhecke}, and the integrand in \eqref{trans211} becomes
\begin{equation} \label{trans212}
\sum_{m<z} \frac{\mu(m)}{m^{2it}} \sum_{l<\frac{z}{m}} \nu (ml) \frac{\lambda(l)}{l^{it}} \sum_{\frac{N}{m^2l} < k \leq \frac{M}{m^2l}} \frac{\lambda(k)}{k^{it}}.
\end{equation}
Similarly for $S_{22}(N)$, it suffices to estimate the following sum.
\begin{equation} \label{trans221}
N^{-\frac{1}{4}} \int_{\frac{T}{2}}^{2T} \left| \sum_{N < n \leq M} \sum_{\substack{ad=n \\ a, d < z}} \nu (d) \lambda(ad) (ad)^{-it} \right| \dif t + N^{\frac{3}{4}+\epsilon}.
\end{equation}
Applying the multiplicative property of Hecke eigenvalues, Lemma~\ref{multhecke}, we arrive at something that is completely analogous to \eqref{trans211}.  The integrand of \eqref{trans221} becomes
\begin{equation} \label{trans222}
\sum_{m<z} \frac{\mu(m)}{m^{2it}} \sum_{l<\frac{z}{m}} \nu (ml) \frac{\lambda(l)}{l^{it}} \sum_{\substack{\frac{N}{m^2l} < k \leq \frac{M}{m^2l} \\ k<\frac{z}{m} }} \frac{\lambda(k)}{k^{it}}.
\end{equation}
Note that the only difference between \eqref{trans212} and \eqref{trans222} is the additional restriction on the range of summation of the inner-most sum of \eqref{trans222}. \newline

Now our proof goes the same as that of Lemma~\ref{S1est}.  We first insert Perron's formula, Lemma~\ref{perronform}, into \eqref{trans212} and \eqref{trans222}, and then move the lines of integration to the critical line and then apply the large sieve and mean value theorems of automorphic $L$-functions.  Lemma~\ref{S2est} follows.
\end{proof}

Combining the lemmas of Sections 4 and 5 and sum up all the dyadic intervals, we have proved Theorem~\ref{sumprimresult}.

\section{Notes}

We could have also applied a result of M. Jutila, Theorem 4.6 in \cite{MJ}.  In brief, Jutila's result gives that
\begin{equation*}
S_1(N) + S_2(N) \ll N^{\frac{5}{6}+\epsilon}.
\end{equation*}
This will not essentially affect our final result of Theorem~\ref{sumprimresult}, but our results of the previous section give better estimates. \newline

Moreover, we can now see that where the obstacle lies in trying to attain the majorant of $N^{\frac{3}{4}}$.  It comes from the want of better means to estimate the bilinear forms in Section 4.  We could, in principle, do the similar thing that we did in Section 5, use the mean value theorems of $L$-function and Dirichlet series to estimate the sums of our interest.  However, the known results do not yield better estimates. \newline

We further note that if one assumes the truth of the Montgomery conjecture \cite{HM3} or the Lindel\"{o}f hypothesis, then we get better estimates for the bilinear forms and obtain the expected majorant implied by \eqref{ilsheur}.  Finally, results similar to Theorem~\ref{sumprimresult} can be obtained for cusp forms for congruence subgroups via analogous means.  The corresponding majorant for the Type I sums weakens as the level increases.

\bibliography{biblio}
\bibliographystyle{amsxport}

\vspace*{.1in}
{\sc\small Dept. Math., U.S. Military Academy, West Point, NY 10996 USA \newline
\indent Email Address: {\tt Liangyi.Zhao@usma.edu}

\end{document}